\documentclass[12pt]{article}  
\usepackage{amssymb}
\usepackage{graphicx}
\usepackage{amsmath}
\usepackage{amsfonts}
\usepackage{hyperref}
\usepackage[utf8]{inputenc}
\usepackage[driver=pdftex]{geometry}
\DeclareUnicodeCharacter{0301}{}
\newcommand{\singlespacing}{\let\CS=\@currsize\renewcommand{\baselinestreatch}{1.0}\tiny\CS}
\newcommand{\doublespacing}{\let\CS=\@currsize\renewcommand{\baselinestreatch}{1.5}\tiny\CS}

 \newtheorem{theorem}{Theorem}[section]
 
 \newtheorem{Proposition}{Proposition}[section]
 \newtheorem{Corollary}{Corollary}[section]

 \newtheorem{rem}{Remark}[section]

 \DeclareMathOperator{\Lie}{\mathcal{L}}

\numberwithin{equation}{section}
\begin{document}
\begin{center}
 {\textbf{Geometric Rigidity of Manifolds Admitting $m$-Modified Conformal Vector Fields}}\\
\vskip.1in
Ramesh Sharma$^1$ and Rahul Poddar$^2$\\
\vskip.1in
$^1$University of New Haven, West Haven, CT 06516, USA. E-mail: rsharma@newhaven.edu  \\
$^1$Sri Sathya Sai Institute of Higher Learning, Prasanthi Nilayam, 515134, AP, INDIA.\\
$^2$N$\bar{a}$land$\bar{a}$ University, Rajgir, 803116, Bihar, INDIA. E-mail: r.poddar@nalandauniv.edu.in   \\
\end{center}

\begin{abstract} 
\noindent
We establish certain integral inequalities for $m$-modified conformal vector fields on compact Riemannian manifolds and derive corresponding triviality and rigidity results. We further obtain integrability conditions for such vector fields and establish an extended Bourguignon-Ezin identity. We then prove that a compact Riemannian manifold admitting a non-trivial closed $m$-modified conformal vector field is isometric to a round sphere under either a constant scalar curvature and positive Ricci curvature assumption or a suitable bound on the Ricci curvature involving the first nonzero eigenvalue of the Laplace operator. As a consequence, we obtain an extension of a classical sphere characterization result of Yano and Nagano to the corresponding result for a non trivial $m$-modified conformal vector field assumed closed. Finally, we consider a complete K\"ahler manifold admitting a closed $m$-modified conformal vector field $V$, and show that $V$ is $m$-modified Killing and around a non-zero point of $V$ the manifold splits as the product of the complex Euclidean plane $\mathbb{C}$ and a complete K\"ahler manifold.
\end{abstract}
 


\section{Introduction}
A Riemannian manifold $(M,g)$ is said to be a Ricci soliton if it satisfies the equation
\begin{equation}\label{1.1}
\frac{1}{2}\Lie_V g + Ric = \lambda g,
\end{equation}
where $\Lie$ denotes the Lie-derivative operator, $V$ is a smooth vector field on $M$, $g$ is the Riemannian metric, $Ric$ is the Ricci tensor of $g$, and $\lambda$ is a real constant. If $V$ is the gradient of a smooth function $f$, i.e., $V=\nabla f$, then $(M,g,f,\lambda)$ is called a gradient Ricci soliton.
\noindent
As a natural generalization of Einstein metrics and gradient Ricci solitons, Catino \cite{catino2012generalized} introduced the notion of a generalized quasi-Einstein manifold. More precisely, an $n$-dimensional Riemannian manifold $(M,g)$ is said to be generalized quasi-Einstein if it satisfies
\begin{equation}\label{1.2}
Ric + Hess\ f - \mu, df\otimes df = \lambda g,
\end{equation}
where $\mu$ is a constant and $\lambda$ is a smooth function on $M$. Catino \cite{catino2012generalized} showed that a complete generalized quasi-Einstein manifold with harmonic Weyl tensor and zero radial Weyl curvature is locally a warped product with an $(n-1)$-dimensional Einstein fiber.
\noindent
In particular, taking

$$ \mu=\frac{1}{m},\qquad 0<m\leq\infty, $$
\noindent
equation \eqref{1.2} becomes
\begin{equation}\label{1.3}
Ric + Hess\ f - \frac{1}{m}df\otimes df =\lambda g.
\end{equation}
The tensor

$$ Ric_f^m=Ric+Hess\ f-\frac{1}{m}df\otimes df $$
\noindent
is referred to as the $m$-Bakry-\'Emery Ricci tensor. When $\lambda$ is constant, equation \eqref{1.3} defines an $m$-quasi-Einstein manifold; see Case, Shu and Wei \cite{case2011rigidity} and He, Petersen and Wylie \cite{he2010classification}.
\noindent
A vector-field analogue of \eqref{1.3} was subsequently considered by Limoncu \cite{limoncu} and Barros and Ribeiro, Jr. \cite{barros2012integral}. More precisely, if $v$ denotes the $1$-form metrically dual to the vector field $V$, they considered
\begin{equation}\label{1.4}
\frac{1}{2}\Lie_V g + Ric - \frac{1}{m}v\otimes v = \lambda g,
\end{equation}
where $\lambda$ is a constant. Moreover, in the limiting case $m=\infty$, equation \eqref{1.4} reduces to the classical Ricci soliton equation \eqref{1.1}.
\noindent
For a Riemannian manifold ($M,g$), it is interesting to know what happens when $g$ is both Einstein and $m$-quasi Einstein. In this case, there exists a vector field $V$ on $M$ such that
\begin{eqnarray*}
\frac{1}{2}\Lie_V g-\frac{1}{m}v\otimes v=c g,
\end{eqnarray*}
where $c$ is a constant. Zhang and Chen \cite{zhang} generalized the above setting by replacing the constant $c$ with a smooth function $\sigma\in C^\infty(M)$, thereby introducing the following definition.\\

\noindent
\textbf{Definition:} A vector field $V$ on a Riemannian manifold ($M,g$) is called an \textit{$m$-modified conformal vector field} if 
\begin{equation}\label{1.5}
\frac{1}{2}\Lie_V g-\frac{1}{m}v\otimes v=\sigma g
\end{equation}
where $m$ is a positive constant and $\sigma \in C^{\infty}(M)$.\\

\noindent
When $\sigma$ is constant and $0<m<\infty$, a vector field $V$ satisfying \eqref{1.5} is called an $m$-modified homothetic vector field; in the special case $\sigma=0$, it is called an $m$-modified Killing vector field. $V$ is said to be trivial when $V=0$. In the limiting case $m=\infty$, equation \eqref{1.5} reduces to: $ \Lie_V g=2\sigma g$, and hence recovers the usual notion of a conformal vector field.
\noindent
Zhang and Chen \cite{zhang} established several rigidity results for $m$-modified conformal vector fields. In particular, they proved that a connected compact Riemannian manifold of dimension at least $2$ admits no $m$-modified conformal vector field of constant length. They further showed in Riemannian case that, for Lie groups admitting left-invariant $m$-modified conformal vector fields, the metric is Einstein and non-trivial $m$-quasi-Einstein. Moreover, such Lie groups are solvable and isometric to a real hyperbolic space.
\noindent
In the compact setting, Poddar and Sharma \cite{Poddar-Sharma} investigated $m$-modified conformal vector fields and established several rigidity results.
In particular, they proved that:
\begin{enumerate}
    \item[(i)] a compact Riemannian manifold admits no non-trivial
    $m$-modified homothetic vector field; and

    \item[(ii)] if $V$ is an $m$-modified conformal vector field, with
    $0<m<\infty$, on a compact Riemannian manifold satisfying
    \[
        \operatorname{Ric}(V,V)
        \leq
        \frac{2(1-n)}{nm^2}
        \big(\max_M |V|\big)^4,
    \]
    then $V$ is trivial.
\end{enumerate}
\noindent
More recently, Deshmukh et al. \cite{Deshmukh1} obtained further triviality
results for $m$-modified conformal vector fields in both compact and
non-compact settings. They also constructed an example of an $m$-modified
conformal vector field on a round sphere.\\

\noindent
The purpose of the present paper is to generalize the results of Poddar and Sharma \cite{Poddar-Sharma} described in (i) and (ii), and to further investigate the geometric and analytic restrictions imposed on an $m$-modified conformal vector field by the compactness of the underlying Riemannian manifold. In particular, we obtain an extension of the classical result of Yano \cite{yano1952harmonic} for conformal vector fields to the setting of $m$-modified conformal vector fields. We also derive several integrability conditions and establish an extension of the Bourguignon-Ezin identity. Further, we prove that a compact Riemannian manifold admitting a non-trivial closed $m$-modified conformal vector field is isometric to a round sphere under either of the following conditions:
the scalar curvature is constant and the Ricci curvature is positive, or the Ricci curvature satisfies a suitable bound involving the first non-zero eigenvalue of the Laplace operator. Finally, we consider a complete K\"ahler manifold admitting a closed $m$-modified conformal vector field $V$, and show that $V$ is $m$-modified Killing and around a non-zero point of $V$ the manifold splits as the product of the complex Euclidean plane $\mathbb{C}$ and a complete K\"ahler manifold.

\section{$m$-Modified Conformal Vector Fields On Compact Riemannian Manifolds}
First, we establish the following generalization of the result (i) of  \cite{Poddar-Sharma} stated earlier.
\begin{theorem}
If $V$ is an $m$-modified conformal vector field on a compact (without boundary) Riemannian manifold ($M,g$), then $\int_{M}(\Lie_V \sigma)dv_g \le 0$, equality holds when $V$ is Killing for $m= \infty$, and $V$ vanishes for $m \ne \infty$.
\end{theorem}
\begin{rem} For $\sigma$ constant, this theorem reduces to the result (i) of  \cite{Poddar-Sharma}.
\end{rem}

\noindent
\textbf{Proof.} First, tracing equation (\ref{1.5}) gives
\begin{equation}\label{2.1}
    div V = \frac{1}{m}|V|^{2} + n\sigma.
\end{equation}
Equation (\ref{1.5}) can also be exhibited as
\begin{equation}\label{2.2}
    g(\nabla_X V,Y) + g(\nabla_Y V, X) = \frac{2}{m}g(V,X)g(V,Y) + 2\sigma g(X,Y).
\end{equation}
The exterior derivative $dv$ of the $1$-form $v$ is given by
\begin{equation}\label{2.3}
    g(\nabla_X V,Y) - g(\nabla_Y V, X) = 2(dv)(X,Y).
\end{equation}
As $dv$ is skew-symmetric, we define a tensor field of type $(1,1)$ given by
\begin{equation}\label{2.4}
    dv(X,Y) = g(FX,Y).
\end{equation}
Evidently, $F$ is skew self-adjoint, i.e., $g(FX,Y)=-g(X,FY)$. Thus, (\ref{2.4}) assumes the form
\begin{equation}\label{2.5}
    g(\nabla_X V,Y) - g(\nabla_Y V, X) = 2g(FX,Y).
\end{equation}
Adding equations (\ref{2.2}) and (\ref{2.5}) side by side, and factoring out $Y$ gives
\begin{equation}\label{2.6}
    \nabla_X V = \frac{1}{m}g(V,X)V + \sigma X + FX.
\end{equation}
Now, we compute the covariant derivative of the squared $g$-norm of $V$ using (\ref{2.6}) as follows
\begin{equation}\label{2.7}
    \nabla_X |V|^{2} = 2g(\nabla_X V, V) = \frac{2}{m}g(V,X)|V|^{2} + 2\sigma g(X,V) + 2g(FX,V).
\end{equation}
Using this we compute
\begin{equation}\label{2.8}
    \Delta |V|^{2} = \frac{2}{m}div(|V|^{2}V) + 2( \Lie_V \sigma +\sigma div V) - 2div (FV).
\end{equation}
Next, multiplying both sides of  (\ref{2.1}) with $div V$, gives $(div V)^2=n\sigma div V +\frac{1}{m}|V|^2 div V$. Integrating it over the compact $M$, using divergence theorem, and the identity: $div(|V|^2 V)=\nabla_V |V|^2 +|V|^2 div V$, we obtain
\begin{equation}\label{2.9}
\int_M (div V)^2 dv_g = \int_M [n\sigma div V-\frac{1}{m}(\nabla_V |V|^2)]dv_g
\end{equation}
where $dv_g$ is the $g$-volume elemment of ($M,g$). Substituting $V$ for $X$ in (\ref{2.7}) we have
\begin{equation}\label{2.10}
\nabla_V |V|^2 =\frac{2}{m}|V|^4 +2\sigma |V|^2 .
\end{equation}
Using this equation, along with the identity: $div(|V|^2 V)=\nabla_V |V|^2 +|V|^2 div V$, in equation (\ref{2.8}) provides
\begin{equation}\label{2.11}
   \frac{1}{2} \Delta |V|^{2} = \frac{1}{m}[\frac{2}{m}|V|^4 +2\sigma |V|^2+|V|^2 div V]+ \Lie_V \sigma +\sigma div V - div (FV).
\end{equation}
Eliminating $|V|^{2}$ between (\ref{2.1}) and (\ref{2.11}) shows that
\begin{equation*}
   \frac{1}{2} \Delta |V|^{2} = \frac{2}{m}|V|^2 (divV-n\sigma)+(2\sigma +div V)(div V-n\sigma)+ \Lie_V \sigma +\sigma div V - div (FV).  
\end{equation*}
Now using (\ref{2.1}) in the above, and simplifying  we find that
\begin{equation*}
   \frac{1}{2} \Delta |V|^{2} = 3(div V)^2 +(\sigma div V)(3-5n)+2n(n-1)\sigma^2 + \Lie_V \sigma - div (FV).  
\end{equation*}
Integrating it over compact $M$ provides
\begin{equation}\label{2.12}
\int_M [3(divV)^2 +(3-5n)(\sigma div V)+2n(n-1)\sigma^2 +\Lie_V \sigma]dv_g=0 .
\end{equation}
Going back to equation (\ref{2.9}), using (\ref{2.10}) and (\ref{2.1}) we have that
\begin{equation}\label{2.13}
\int_M [3(divV)^2 +(2-5n)(\sigma div V)-2n(1-n)\sigma^2]dv_g=0.
\end{equation}
At this point, multiplying the above equation with $3-5n$, and (\ref{2.12}) with $2-5n$, and subsequently subtracting the resulting equations we obtain
\begin{equation}\label{2.14}
\int_M [3(divV)^2 +2n(n-1)\sigma^2]dv_g=(2-5n)\int_M (\Lie_V \sigma)dv_g .
\end{equation}
The assertions of the theorem follow from equation (\ref{2.14}). This completes the proof.\\

\noindent
We turn our attention back to the result (ii) of \cite{Poddar-Sharma} that provides a generalization of the following result \cite{yano1952harmonic} of Yano ``If a compact Riemannian manifold $(M,g)$ admits a conformal vector field $V$ satisfying $Ric(V,V) \leq 0$, then $V$ is parallel, and if the Ricci curvature is negative definite, then the conformal vector field $V$ is zero". Improving this result further and extending the aforementioned result of Yano to the case of an $m$-modified vector field, we prove the following result.

\begin{theorem}
Let ($M,g$) be a compact Riemannian manifold admitting an $m$-modified conformal vector field $V$. If $Ric(V,V) \le 0$, then $V$ is parallel if $m =\infty$, and $V$ vanishes identically if $m \ne \infty$.
\end{theorem}

\noindent
\textbf{Proof.} Recalling (\ref{2.6}) and using it we compute
\begin{eqnarray}\label{2.15}
Riem(Y,X)V&=&(Y \sigma)X-(X \sigma)Y+\frac{1}{m}[g(\nabla_Y V,X)V-g(\nabla_X V,Y)V\nonumber\\
&+&\sigma v(X)Y-\sigma v(Y)X+v(X)FY\nonumber\\
&-&v(Y)FX]+(\nabla_Y F)X-(\nabla_X F)Y,
\end{eqnarray}
where $Riem$ denotes the curvature tensor of type (1,3) for $g$, defined by: $Riem(X,Y)Z=\nabla_X \nabla_Y Z-\nabla_Y \nabla_X Z-\nabla_{[X,Y]}Z$ for arbitrary $X,Y,Z \in \mathfrak{X}(M)$.
Contracting equation (\ref{2.15}) at $Y$ we have
\begin{eqnarray*}
Ric(X,V)&=&(1-n)X\sigma +\frac{1}{m}[g(\nabla_V V,X)-g(\nabla_X V,V)\nonumber\\
&+&(n-1)\sigma v(X)-g(FX,V)]+(div F)X.
\end{eqnarray*}
Next, substituting $V$ for $X$ in the above equation gives
\begin{equation}\label{2.16}
Ric(V,V)=(1-n)V\sigma +\frac{1}{m}(n-1)\sigma|V|^2 +(divF)V
\end{equation}
Through the computation: 
\begin{eqnarray*}
\nabla_X(FV)&=&(\nabla_X F)V+F\nabla_X V\nonumber\\
&=& (\nabla_X F)V+\sigma FX +\frac{1}{m}v(X)FV+F^2 X
\end{eqnarray*} and its trace at $X$ we obtain
\begin{equation*}
div(FV)=(div F)V -|F|^2,
\end{equation*}
where we have used (\ref{2.6}) and the skew self adjointness of $F$. Integrating the preceding equation over the compact $M$, we have
\begin{equation}\label{2.17}
\int_M (div F)V dv_g=-\int_M |F|^2 dv_g
\end{equation}
Integrating equation (\ref{2.16}) over the compact $M$ and making use of equations (\ref{2.17}) and (\ref{2.14}) we obtain the formula
\begin{eqnarray}\label{2.18}
\int_M Ric(V,V)dv_g &=& \frac{n-1}{5n-2}\int_M [3(div V)^2 +2n(n-1)\sigma^2]dv_g\nonumber\\
&+&\frac{n-1}{m}\int_M |V|^2 dv_g +\int_M |F|^2dv_g.
\end{eqnarray}
The conclusions of the theorem follow immediately from this consequence, completing the proof.

\section{Integrability Conditions for $m$-Modified Conformal Vector Field} In this section, we derive integrability conditions for an $m$-modified conformal vector field $V$ on a Riemannian manifold ($M,g$). By definition, the equation (\ref{1.5}) is written in local co-ordinates as
\begin{equation}\label{3.1}
\Lie_V g_{ij}=2\sigma g_{ij}+\frac{2}{m}V_i V_j.
\end{equation}
Lie-differentiating the relation:  $g_{ij}g^{jk}=\delta^k_i$ along $V$ gives
\begin{equation*}
(\Lie_V g_{ij})g^{jk}+g_{ij} \Lie_V g^{jk}=0.
\end{equation*}
Using (\ref{3.1}) in it, and taking inner product with $g^{il}$ yields that
\begin{equation*}
\Lie_V g^{kl}=-2g^{kl}-\frac{2}{m}V^k V^l .
\end{equation*}
At this point we denote $V_i V_j$ by $p_{ij}$ and hence $Tr. p=|V|^2$. Then we have
\begin{equation}\label{3.2}
\Lie_V g_{ij}=2\sigma g_{ij}+\frac{2}{m}p_{ij}, \hskip.1 in \Lie_V g^{ij}=-2\sigma g^{ij}-\frac{2}{m}p^{ij}.
\end{equation}
Use of the formula (Yano \cite{yano1970integral}):
\begin{equation*}
\Lie_V \Gamma^i _{jk}=\frac{1}{2}g^{il}[\nabla_j \Lie_V g_{kl}+\nabla_k \Lie_V g_{jl}-\nabla_l \Lie_V g_{jk}]
\end{equation*}
in (\ref{3.2}) provides
\begin{eqnarray}\label{3.3}
\Lie_V \Gamma^i _{jk}&=&(\nabla_j \sigma)\delta^i _k +(\nabla_k \sigma)\delta^i _j -(\nabla^i \sigma)g_{jk}\nonumber\\
&+&\frac{1}{m}[\nabla_j  p^i _k +\nabla_k  p^i _j -\nabla^i  p_{jk}]
\end{eqnarray}
where $\Gamma^i _{jk}$ denote the connection coefficients of $g_{ij}$. Now combined use of the formula \cite{yano1970integral}: $\Lie_V R_{ikl}\ ^j=\nabla_i \Lie_V \Gamma^j _{kl}-\nabla_k \Lie_V \Gamma^j _{il}$ and equation (\ref{3.3}), gives
\begin{eqnarray}\label{3.4}
\Lie_V R_{ikl}\ ^j&=&(\nabla_i \nabla_l \sigma)\delta^j _k -(\nabla_k \nabla_l \sigma)\delta^j _i-(\nabla_i \nabla^{j}\sigma)g_{kl}+(\nabla_k \nabla^{j}\sigma)g_{il}\nonumber\\ 
 &+&\frac{1}{m}[\nabla_i \nabla_k p^j _l -\nabla_k \nabla_i p^j _l+\nabla_i \nabla_l p^j _k\nonumber\\
&-&\nabla_k \nabla_l p^j _i -\nabla_i \nabla^{j}p_{kl}+\nabla_k \nabla^{j}p_{il}]
\end{eqnarray}
where $R_{ikl}\ ^j X^l=\nabla_i \nabla_k X^j -\nabla_k \nabla_i X^j$ and $X^l$ are components of an arbitrary vector field $X \in \mathfrak{X}(M)$. Contracting (\ref{3.4} with respect to $i$ and $j$ we get
\begin{eqnarray}\label{3.5}
\Lie_V R_{kl}&=&(2-n)\nabla_k \nabla_l \sigma -(\Delta \sigma)g_{kl}\nonumber\\
&+&\frac{1}{m}[\nabla_i \nabla_k p^i _l+\nabla_i \nabla_l p^i _k -\nabla_k \nabla_l (Tr. p)]
\end{eqnarray}
where $R_{kl}$ are the components of the Ricci tensor of $g_{ij}$. Finally, we compute $\Lie_V R=\Lie_V(R_{kl}g^{kl})$ using (\ref{3.2}) and (\ref{3.5}) and thus obtain
\begin{eqnarray}\label{3.6}
\Lie_V R&=&2(1-n)\Delta \sigma -2\sigma R+\frac{2}{m}[\nabla_i \nabla^{l}p^i _l\nonumber\\
&-&\Delta (Tr.p)-R_{kl}p^{kl}].
\end{eqnarray}
Equations (\ref{3.3}),(\ref{3.4}),(\ref{3.5}),(\ref{3.6}) provide the evolutions of the connection coefficients, curvature tensor, Ricci tensor and scalar curvature along an $m$-modified conformal vector field $V$.
\section{Extended Bourguignon-Ezin Identity and Scalar Curvature Constraints}
Here we assume ($M,g$) compact and admitting an $m$-modified vector field $V$. Integrating (\ref{3.6}) over the compact $M$ we find that
\begin{equation}\label{4.1}
\int_M (\Lie_V R)dv_g=-2\int_M \bigg[\sigma R+\frac{1}{m}Ric(V,V)\bigg]dv_g
\end{equation}
Using the equation (\ref{2.1}) and the identity: $div(RV)=\Lie_V R+Rdiv(V)$ in (\ref{4.1}) we obtain 
\begin{equation}\label{4.2}
\int_M (\Lie_V R)dv_g=-\frac{2n}{m(n-2)}\int_M\bigg[\bigg(Ric -\frac{R}{n}g\bigg)(V,V)\bigg]dv_g.
\end{equation}
We note that this reduces to the Bourguignon-Ezin identity \cite{Bourguignon-Ezin} for a conformal vector field (for which $m =\infty$) on a compact Riemannian manifold. We state the above consequence as the following result.
\begin{theorem}For an $m$-modified conformal vector field $V$ on a compact Riemannian manifold of dimension $n > 2$, we have the extended Bourguignon-Ezin identity given by (\ref{4.2}).
\end{theorem}

\begin{rem} Equation (\ref{4.2}) is trivially satisfied by a compact Einstein manifold. 
\end{rem}

\noindent
We close this section with the scalar curvature constraint. Using equation (\ref{2.18}) in (\ref{4.2}) and simplifying gives the following constraint on the scalar curvature $R$:
\begin{eqnarray}\label{4.3}
&&-\frac{m(n-2)}{2n}\int_M (\Lie_V R)dv_g =\frac{n-1}{5n-2}\int_M [3(div V)^2 +2n(n-1)\sigma^2]dv_g\nonumber\\
&+&\int_M |F|^2 dv_g +\int_M \frac{n-1}{m}|V|^2 dv_g-\frac{1}{n}\int_M R|V|^2 dv_g.
\end{eqnarray}
It follows immediately from the above equation that
\begin{equation}\label{4.4}
\int_M (\Lie_V R)dv_g \le \frac{2}{m(n-2)}\int_M \bigg(R-\frac{n(n-1)}{m}\bigg)|V|^2.
\end{equation}
In particular, if $R$ is constant, and $V$ is non-zero and non-conformal, then the inequality (\ref{4.4}) shows that $R \ge \frac{n(n-1)}{m}$, because $V$ is not conformal, i.e., $m \ne \infty$. However, if $R=\frac{n(n-1)}{m}$, then (\ref{4.3}) implies that $div V=0, \sigma =0$, and hence from (\ref{2.1}) we get $V=0$, a contradiction. Consequently,  $R > \frac{n(n-1)}{m}$. Hence we have proved the following result.
\begin{Proposition} If a compact Riemannian manifold of dimension $n>2$ has constant scalar curvature $R$, and admits a non-zero non-conformal $m$-modified conformal vector field, then $R$ is bounded below by the positive constant  $\frac{n(n-1)}{m}$.
\end{Proposition}

\section{Closed $m$-modified Conformal Vector Fields On Compact Riemannian Manifolds}

We begin this section by recalling the following rigidity theorems for conformal vector fields on compact Riemannian manifolds.
\begin{theorem}[Yano \cite{yano1970integral}]
    Let $(M,g)$ be an $(n\geq 3)$-dimensional closed Riemannian manifold with constant scalar curvature admitting a nontrivial gradient conformal vector field, then $M$ is isometric to a round sphere.
\end{theorem}
\begin{theorem}[Deshmukh and Al-Solamy \cite{Deshmukh}]
Let $(M,g)$ be an $n$-dimensional compact Riemannian manifold whose Ricci curvature satisfies
\[
0 < \operatorname{Ric}
\leq
(n-1)\left(2-\frac{nc}{\lambda_1}\right)c,
\]
where $c$ is a constant and $\lambda_1$ denotes the first nonzero eigenvalue of the Laplace operator. If $M$ admits a nonzero conformal
gradient vector field, then $M$ is isometric to the sphere $\mathbb{S}^n(c)$.
\end{theorem}

\noindent
Motivated by the aforementioned results, we now establish the corresponding rigidity result (a generalization of Theorem 5.2) for closed $m$-modified conformal vector fields.

\begin{theorem}\label{thm:closed-m-modified}
Let $(M^n,g)$ be a compact Riemannian manifold with positive Ricci curvature admitting a non-trivial closed $m$-modified conformal vector
field. Then $(M,g)$ is isometric to a round sphere if one of the following conditions is satisfied:
\begin{enumerate}
    \item[(i)] $M$ has constant scalar curvature and $n\geq3$.
    
    \item[(ii)] The Ricci curvature satisfies
    \[
    \operatorname{Ric}
    \leq
    (n-1)\left(2-\frac{nc}{\lambda_1}\right)c,
    \]
    where $c$ is a constant and $\lambda_1$ denotes the first nonzero
    eigenvalue of the Laplace operator.
\end{enumerate}
\end{theorem}
\noindent
\textbf{Proof.}
Let $V$ be a non-trivial closed $m$-modified conformal vector field on $M$, and let $v=V^\flat$ denote its metrically dual $1$-form. Since $M$ is compact, the Hodge decomposition gives
\begin{equation}\label{eq:hodge}
    v=dh+\eta,
\end{equation}
where $h\in C^\infty(M)$ and $\eta$ is a $1$-form such that $\delta \eta = 0$. We are going to show that the harmonic component in \eqref{eq:hodge} vanishes under our hypothesis. Since $V$ is closed, its dual $1$-form $v$ is closed, and hence
\[
    0=dv=d\eta.
\]
On the other hand, the Hodge decomposition gives
\[
    \delta\eta=0.
\]
Thus
\[
    \Delta_H\eta=(d\delta+\delta d)\eta=0,
\]
so that $\eta$ is harmonic. The Bochner formula for a harmonic $1$-form yields
\begin{equation}\label{eq:bochner-harmonic}
    \frac12\Delta |\eta|^2
    =
    |\nabla\eta|^2+\operatorname{Ric}(\eta^\sharp,\eta^\sharp),
\end{equation}
where $\eta^\sharp$ denotes the vector field metrically dual to $\eta$. Integrating \eqref{eq:bochner-harmonic} over the compact manifold $M$ and applying the divergence theorem, we obtain
\begin{equation}\label{eq:integrated-bochner}
    0
    =
    \int_M
    \left(
        |\nabla\eta|^2+
        \operatorname{Ric}(\eta^\sharp,\eta^\sharp)
    \right)\,d\mu.
\end{equation}
Since $\operatorname{Ric}$ is positive definite, both terms in the integrand are non-negative and
\[
    \operatorname{Ric}(\eta^\sharp,\eta^\sharp)>0
\]
whenever $\eta^\sharp\neq0$. Consequently, \eqref{eq:integrated-bochner} implies
\[
    \eta\equiv0.
\]
It follows from \eqref{eq:hodge} that
\[
    v=dh,
\]
and hence
\begin{equation}\label{eq:V-gradient}
    V=\nabla h.
\end{equation}
\noindent
We next exploit the $m$-modified conformal equation. As $V=\nabla h$, equation (\ref{1.5}) becomes
\begin{equation}\label{eq:h-equation}
    \operatorname{Hess}h
    -\frac1m dh\otimes dh
    =\sigma g.
\end{equation}
Let us define: $ u=e^{-h/m}$. Then
\[
    du=-\frac1m e^{-h/m}dh
\]
and a direct computation gives
\begin{equation}\label{eq:hess-u-computation}
    \operatorname{Hess}u
    =
    \frac{e^{-h/m}}{m^2}\,dh\otimes dh
    -\frac{e^{-h/m}}{m}\operatorname{Hess}h.
\end{equation}
Using \eqref{eq:h-equation} in \eqref{eq:hess-u-computation}, we obtain
\[
\begin{aligned}
    \operatorname{Hess}u
    &=
    \frac{e^{-h/m}}{m^2}\,dh\otimes dh
    -\frac{e^{-h/m}}{m}
      \left(
        \frac1m dh\otimes dh+\sigma g
      \right)\\
    &=
    -\frac1m e^{-h/m}\sigma g.
\end{aligned}
\]
Therefore,
\begin{equation}\label{eq:hess-u-lambda}
    \operatorname{Hess}u=\lambda g,
    \qquad
    \lambda=-\frac1m e^{-h/m}\sigma.
\end{equation}
Thus $\nabla u$ is a gradient conformal vector field. It remains to verify that $\nabla u$ is non-trivial. Suppose, to the contrary, that $\sigma = 0$. Then integrating equation (\ref{2.1}) and using the divergence theorem would imply $V\equiv 0$, a contradiction to the nontriviality of $V$. Therefore, $\sigma$ is a non-zero function on $M$. Hence $(M,g)$ admits a non-trivial gradient conformal vector field. To prove part (i), we assume that $M$ has constant scalar curvature. Since $n\geq 3$, Theorem 5.1 applies and yields that $(M,g)$ is isometric to a round sphere. This completes the proof of part (i).\\
\noindent
For part (ii), the positivity of the Ricci curvature, together with the assumed upper bound on the Ricci curvature, allows us to apply Theorem 5.2. Consequently, $(M,g)$ is isometric to a round sphere. This completes the proof.
\noindent
\begin{Corollary}
Let $(M^n,g)$ be an $n$-dimensional compact Riemannian manifold, whose Ricci curvature satisfies
\[
0 < \operatorname{Ric}
\leq
(n-1)\left(2-\frac{nc}{\lambda_1}\right)c,
\]
for a constant $c$, where $\lambda_1$ is the first nonzero eigenvalue of the Laplace operator. If $M$ admits a non-trivial closed conformal vector field, then $M$ is isometric to a round sphere.
\end{Corollary}
\noindent
\textbf{Proof.} The result follows immediately by letting $m=\infty$ in (\ref{1.5}) and applying Theorem 5.3 (ii).
\begin{Corollary}
Let $(M^n,g)$, $n\geq3$, be a compact Einstein manifold admitting a non-trivial closed $m$-modified conformal vector field ($m \ne \infty$). Then $(M,g)$ is isometric to a round sphere.
\end{Corollary}
\noindent
\textbf{Proof.}
The result follows directly from Theorem 5.3 (i) and Proposition 4.1 (here non-conformal means ($m \ne \infty$)), together with the fact that an Einstein manifold has constant scalar curvature for $n>2$, as a consequence of Schur's lemma.
\noindent
\begin{rem}
    Corollary 5.2 raises the following question ``What can we say about a compact Einstein manifold admitting a non-trivial non-closed $m$-modified conformal vector field ($m \ne \infty$)?" For $m =\infty$ ($V$ conformal) we know, by a classical result of Yano and Nagano \cite{Yano-Nagano}, that a compact Einstein manifold $M$ admitting a non-Killing conformal vector field $V$ is isometric to a Euclidean sphere (this result was proved more generally in \cite{Yano-Nagano}, for a complete $M$ under the condition that $V$ is complete. For $M$ compact, $V$ is automatically complete.)
\end{rem}

\section{Closed $m$-modified Conformal Vector Fields On K\"ahler Manifolds} An almost Hermitian manifold $M$ is a real $2n$-dimensional smooth manifold with a $(1,1)$ tensor field $J$ and a Riemannian metric $g$ such that  $J{^2} = -I$ and $g(JX,JY) = g(X,Y)$ for arbitrary vector fields $X$, $Y$ on $M$. If $J$ is integrable, i.e., the Nijenhuis tensor $N$ defined by $N[X,Y] = [JX, JY] - [X,Y] - J[JX, Y] - J[X, JY]$ vanishes, then $M$ is a Hermitian manifold. An almost Hermitian manifold is said to be an almost  K\"ahler manifold if the fundamental $2$-form $\Omega$ defined by $\Omega(X,Y)= g(JX,Y)$ is closed i.e., $d\Omega = 0$ ($d$ is the exterior derivative operator). An almost  K\"ahler manifold whose underlying almost complex structure $J$ is integrable i.e., the Nijenhuis tensor of $J$ vanishes on $M$, is known as a  K\"ahler manifold. It is known that for a  K\"ahler manifold $\nabla J = 0$, where $\nabla$ is the Riemannian connection of $g$. Conformal vector fields on  K\"ahler manifolds are known to exhibit strong rigidity phenomena. A classical result of Goldberg \cite{Goldberg} says that a closed conformal vector field (i.e., the metrically equivalent $1$-form is closed) on a  K\"ahler manifold is necessarily homothetic.  In this section, we examine the impact of the existence of a non-zero closed $m$-modified conformal vector field $V$ on a K\"ahler manifold, and prove the following result.
\begin{theorem}If a complete K\"ahler manifold $(M,g,J)$ of real dimension $2n > 2$ admits a real non-trivial $m$-modified ($m \ne \infty$) closed conformal vector field $V$, then $V$ reduces to $m$-modified ($m \ne \infty$) Killing vector field, and around any non-zero of $V$, $M$ splits as the product $\mathbb{C}\times N$, where $\mathbb{C}$ is the Euclidean complex plane and $N$ is a complete K\"ahler manifold of real dimension $2n-2$.
\end{theorem}
\textbf{Proof.} The $m$-modified equation can be written in view (\ref{2.6}) and (\ref{2.4}) as
\begin{equation}\label{6.1}
\nabla_X V=\frac{1}{m}g(V,X)V+\sigma X
\end{equation}
for arbitrary $X \in \mathfrak{X}(M)$, because $V$ is closed, i.e. $dv=0$ and hence $F=0$. Operating (\ref{6.1}) by the complex structure tensor $J$ we have
\begin{equation}\label{6.2}
J\nabla_X V=\frac{1}{m}g(V,X)JV+\sigma JX
\end{equation}
Next, substituting $JX$ for $X$ in (\ref{6.1}) we also have
\begin{equation}\label{6.3}
\nabla_{JX}V=\frac{1}{m}g(V,JX)JV+\sigma JX
\end{equation}
Before proceeding ahead, we denote the $1$-form metrically equivalent to $JV$ by $(JV)^\flat$, i.e., $(JV)^{\flat}(X)=g(JV,X)$. Using (\ref{6.2}) and (\ref{6.3}) we find that
\begin{eqnarray*}
2d(JV)^{\flat}(X,Y)&=&g(\nabla_X JV,Y)-g(\nabla_Y JV,X)\nonumber\\
&=&g(J\nabla_X V,Y)-g(J\nabla_Y V,X)
\end{eqnarray*}
The use of equation (\ref{6.2}) in the above equation entails
\begin{equation}\label{6.4}
2(d (JV)^{\flat})(X,Y)=\frac{1}{m}[g(V,X)g(JV,Y)-g(JV,X)g(V,Y)]+2\sigma g(JX,Y),
\end{equation}
where $d$ is the exterior derivative operator and $X, Y$ arbitrary vector fields on $M$. Skew-symmetrizing it with respect to $X,Y$ shows that
\begin{equation}\label{6.5}
2 d(JV)^{\flat}=\frac{1}{m}V^{\flat} \wedge (JV)^{\flat}+2\sigma \Omega.
\end{equation}
Operating (\ref{6.5}) with $d$ gives
\begin{equation}\label{6.6}
\frac{1}{m}V^{\flat} \wedge d (JV)^{\flat}=2 d\sigma \wedge \Omega
\end{equation}
where we used the co-boundary formula: $d^2=0$, $dv=0$, and the property: $d(\alpha \wedge \beta)=d\alpha \wedge \beta +(-1)^p \alpha \wedge d\beta$ for $p$ and $q$ forms $\alpha, \beta$ respectively. Expanding (\ref{6.6}) for arbitrary $X,Y,Z \in \mathfrak{X}(M)$ we have
\begin{eqnarray*}
&&v(X)d (JV)^{\flat}(Y,Z)+v(Y)d(JV)^{\flat}(Z,X)+v(Z)d(JV)^{\flat}(X,Y)\nonumber\\
&=&2m[(X\sigma) \Omega(Y,Z)+(Y\sigma) \Omega(Z,X)+ (Z\sigma) \Omega(X,Y)].
\end{eqnarray*}
Using (\ref{6.4}) in the foregoing equation and factoring out $Z$ provides
\begin{equation*}
\sigma[g(V,X)JY-g(V,Y)JX+g(JX,Y)V]=2m[(X\sigma)JY-(Y\sigma)JX+g(JX,Y)\nabla \sigma].
\end{equation*}
Contracting it with respect to $Y$ we get
\begin{equation}\label{6.7}
\nabla \sigma=\frac{\sigma}{2m}V
\end{equation}
because $2n > 2$, by hypothesis. Now, using (\ref{6.1}) we compute and find that
\begin{equation}\label{6.8}
Riem(Y,X)V=(Y\sigma)X-(X\sigma)Y+\frac{\sigma}{m}[v(X)Y-v(Y)X].
\end{equation}
Contracting it at $Y$ and using (\ref{6.7}) we find that
\begin{equation*}
Ric(V)=(2n-1)\nabla \sigma.
\end{equation*}
where we have denoted the Ricci tensor of type (1,1) also by $Ric$. As $Ric$ commutes with $J$ on a K\"ahler manifold, the foregoing equation implies that
\begin{equation}\label{6.9}
Ric (JV)=(2n-1)J\nabla \sigma.
\end{equation}
Going back to (\ref{6.8}), operating by $J$ and using the K\"ahler property: $JRiem(X,Y)Z=Riem(X,Y)JZ$ for all $X,Y,Z \in \mathfrak{X}(M)$, and then contracting at $Y$ gives
\begin{equation}\label{6.10}
Ric(JV)=-J\nabla \sigma +\frac{\sigma}{m}JV
\end{equation}
Comparing equations (\ref{6.9}) and (\ref{6.10}) and noting that $2n>2$ and $m \ne \infty$ by hypothesis, shows that
\begin{equation}\label{6.11}
\sigma V=0.
\end{equation}
Using this in (\ref{6.7}) at once implies that $\sigma$ is constant. Therefore, as $V \ne 0$, (\ref{6.11}) implies that $\sigma=0$. Hence (\ref{6.1}) reduces to 
\begin{equation}\label{6.12}
\nabla_X V=\frac{1}{m}g(V,X)V.
\end{equation}
Now, as $g$ is Hermitian, $g(JX,JY)=g(X,Y)$, we have $g(JX,Y)=-g(X,JY)$ for arbitrary vector fields $X,Y \in \mathfrak{X}(M)$. So $g(JV,V)=-g(V,JV)=0$, and hence $JV$ is orthogonal to $V$. On an open dense subset $\mathcal{O}$ of $M$, let $V \ne 0$, and $(V,JV)$  denote the $2$-dimensional distribution spanned by $V$ and $JV$. To show that $(V,JV)$ is integrable, it suffices to show that $[V,JV] \in (V,JV)$. For this, $[V,JV]=\nabla_V JV-\nabla_{JV}V=J\nabla_V V-\nabla_{JV}V=\frac{1}{m}|V|^2 JV$, where we used (\ref{6.12}).
Now, for real functions $a,b$ on $\mathcal{O}$, computing
\begin{eqnarray*}
\nabla_{aV+bJV}(cV+dJV)&=& a[(Vc)V+\frac{c}{m}|V|^2 V+(Vd)JV \nonumber\\
&+&\frac{d}{m}|V|^2]+b[(JV)c V+(JV)d (JV)]
\end{eqnarray*}
shows that $(V,JV)$ is totally geodesic in $\mathcal{O}$. Next we consider the distribution $(V,JV)^{\perp}$ orthogonal to $(V,JV)$ and complementary in the tangent bundle on $\mathcal{O}$. Let $X,Y \in (V,JV)^{\perp}$, i.e., $g(X,W)=0$, $g(Y,W)=0$ for any $W \in (V,JV)$, i.e. $W=aV+bJV$, for arbitrary real functions $a,b$ on $\mathcal{O}$. Then we have $g(\nabla_Y X,W)+g(X,\nabla_Y W)=0$ and $g(\nabla_X Y,W)+g(Y,\nabla_X W)=0$. Subtracting them and using (\ref{6.12}), we get $g([X,Y],W)=0$. Thus $[X,Y] \in (V,JV)^{\perp}$. Finally, to show that $(V,JV)^{\perp}$ is totally geodesic in $\mathcal{O}$, we compute
\begin{eqnarray*}
g(\nabla_X Y,aV+bJV)&=&-ag(\nabla_X V,Y)-bg(J\nabla_X V,Y)\nonumber\\ 
&-&\frac{a}{m}v(X)v(Y)-\frac{b}{m}v(X)g(JV,Y)=0,
\end{eqnarray*}
because $X,Y$ are $\perp$ to $V,JV$. Also, from (\ref{6.8}) and the consequence $\sigma =0$, we have $Riem(X,Y)V=0$ [and hence $Riem(X,Y)JV=J Riem(X,Y)V=0$] for any $X,Y \in \mathfrak{X}(M)$. Thus the leaves of the foliation $(V,JV$) are flat. We also note that a totally geodesic submanifold in K\"ahler $M$ of codimension $2$ is also K\"ahler. Summing up, we have shown that around any non-zero of $V$, $M$ splits locally as $\mathbb{C}\times N$ where $N$ is a complete $(2n-2)$-dimensional K\"ahler manifold. This completes the proof.


\section{Statements and Declarations}

\textbf{Competing Interests}: On behalf of all the authors, the corresponding author states that there is no conflict of interest. \\

\noindent
\textbf{Data Availability}: Data sharing is not applicable to this article as no datasets were generated or analysed during the current study.\\

\noindent
\textbf{No Use Of Artificial Intelligence}: No Artificial Intelligence has been used in the preparation of this article.\\

\end{document}